\documentclass[a4paper, 12pt]{article}
\usepackage{amsmath,amsthm, amssymb,color}
\usepackage{multicol}
\usepackage{wrapfig}
\usepackage{graphicx}
\usepackage{graphics}
\usepackage{yfonts}
\usepackage{times}
\usepackage{enumerate}
\usepackage{booktabs}
\usepackage[justification=centering]{caption}
\usepackage[font=footnotesize]{caption}

\numberwithin{equation}{section}

\newtheorem{theorem}{Theorem}[section]

\newtheorem{remark}{Remark}

\DeclareMathOperator{\sgn}{sign}

\def\wt{\widetilde}

\def\ben{\begin{eqnarray}}
\def\een{\end{eqnarray}}
\def\be{\begin{eqnarray*}}
\def\ee{\end{eqnarray*}}
\def\Rn{{\mathbb R}^n}
\newcommand{\dive}{\operatorname{div}}
\def\mean#1{{\left\{\!\left|#1\right|\!\right\}}}
\def\cP{{\mathcal P}}
\def\dvg{{\rm div}}

\begin{document}

\title{Functional a posteriori error estimates\\
for parabolic obstacle problems}

\author{{\textsc{ Darya E.  Apushkinskaya}}\\
Department of Mathematics, Saarland University, P.O. Box 151150\\
66041 Saarbr{\"u}cken, Germany\\
\textit{E-mail: darya@math.uni-sb.de}\\
\\
{\textsc{ Sergey I. Repin}}\\
Steklov Institute of Mathematics at St. Petersburg,\\ Fontanka 27, 191023, St. Petersburg,
Russia\\
\textit{E-mail: repin@pdmi.ras.ru}\\
University of Jyv{\"a}skyl{\"a}, P.O. Box 35  (Agora),\\
FIN-40014, Finland \\
\textit{E-mail:serepin@jyu.fi}}

\maketitle

\abstract{The paper is concerned with functional type a posteriori
estimates for the initial boundary value problem for a parabolic partial
differential equation with an obstacle. We deduce a guaranteed and 
computable bound of the distance between the exact minimizer and any
function from the admissible (energy) class of functions. Applications
to the analysis of modeling errors caused by data implification are discussed. An important case of time incremental approximations is specially studied. Numerical examples presented
in the last section show how the estimates work in practice.}


\section{Introduction} 

Mathematical models based on parabolic type equations with obstacles
arise in various branches of science and technology: e.g., 
in mathematical biology (\cite{AABBK2011}), phase transitions problems (\cite{R1971}, \cite{V1996}), electrochemical industry (\cite{El1980}), stochastic control theory (\cite{BL1982}), mathematical economy (\cite{vM1974}, \cite{vM1975}, \cite{PS2007}).
Obstacle problems for parabolic equations are well studied from the mathematical
point of view. Existence of  a generalised solution for the case of a smooth obstacle has been studied in many publications. For time-independent obstacles, first results were obtained in  \cite{LS1967} and \cite{B1972}. The case, where obstacles are presented by  the functions
``regular''  with respect to time was studied in \cite{B1972a}. The case obstacles non--increasing  in time obstacles has been considered in the books \cite{L1969} and \cite{Na1984}. The existence results for linear parabolic problems with general obstacles depending on time only as measurable functions can be found in \cite{MiPu1977}. For irregular obstacles, the comprehensive existence theory was developed in \cite{BDM2011} and \cite{S2015}.
Qualitative properties of solutions and free boundaries for the smooth obstacle case were studied in \cite{Fr1975} and \cite{BlaDoMo2006} in one dimension, and in \cite{C1977}, \cite{ASU2000}-\nocite{ASU2002}\cite{ASU2003}, \cite{CPS2004}, \cite{Bla2006},  \cite{LiMo2015} in higher dimensions, respectively. A systematic overview of the regularity results for smooth obstacles can be found in the book \cite{A2018}. When the obstacle is non-smooth, the regularity properties of solutions and free boundaries were examined in \cite{PS2007}. The regularity of solutions and free boundaries in the so-called parabolic thin obstacle problem (known also as parabolic Signorini problem) was studied in \cite{ArU1988} and \cite{ArU1996} (see also recent publications \cite{DGPT2017}, \cite{BSVGZ2017}, and \cite{Sh2020}).

There exist various numerical methods for solving this class of nonlinear
problems. At this point, we refer to the monographs \cite{G2008}, \cite{T2006}
and the literature cited therein.

Investigation of a priori error estimates for problems with obstacles
begins with the paper by R. Falk \cite{F1974}
devoted to the elliptic case.  Estimates of this type for evolutionary variational
inequalities have been later studied in many papers (e.g., see
\cite{Fe1987} and \cite{V1990}).

In this paper, we discuss a different problem. Our analysis is focused not on properties of the exact minimizer, but on guaranteed  bounds of the difference between the exact solution (minimizer) of the parabolic variational problem and any function (approximation) from the energy class satisfying the prescribed boundary conditions and the restrictions stipulated by the obstacle. They can be called
{\em estimates of deviations} from the exact solution (or
a posteriori estimates of the functional type). The estimates bound
a certain measure (norm) of the error by a functional (error majorant) that depends
on the problem data and approximation type, but do not explicitly depend
on the exact solution. Hence the functional is fully computable and can be used to evaluate the accuracy of an approximation. Within the framework
of this conception, the estimates should be derived on the functional
level by the same tools as commonly used in the theory of partial
differential equations. They do not use specific features of  approximations (e.g., Galerkin orthogonality) what is typical for
a posteriori methods applied in mesh adaptive computations based
upon finite element technologies.
Unlike the a priori rate convergence estimates that establish
general asymptotic properties of an approximation method, these a posteriori  estimates are applied to a particular solution and
allow us to directly verify its accuracy.
 For various elliptic
and parabolic problems estimates of this type have been derived
in \cite{Re2000,Re2002,Re2007} and many subsequent publications. The
reader can find a consequent exposition of the corresponding theory
in the monographs \cite{Re2008} and \cite{RS2020}.
In this paper, we derive such type estimates for the parabolic
obstacle problem. They depend 
only on the approximation solution (which is known) and on the data of the problem. We emphasise that they also
do not need knowledge on the r exact coincidence set associated with
the exact solution. 
The obtained error majorant is non-negative and vanishes if and only if the approximation coincides with the exact minimizer.
It provides a guaranteed bound of the
error expressed in terms of a natural measure of the distance between the
exact and approximate solution for a finite time interval $[0,T]$.

The outline of the paper is as follows. The first part of Section 2 contains basic notation
and the mathematical formulation of the problem. The second part presents main result  (Theorem~1) and discusses it.  In Section~3 we discuss
some applications of the error majorant. First, we show that
it yields simple bounds for modeling errors generated by 
simplification of problem data. The corresponding estimate is
directly computable and do not require an information about the exact solution of the original (complicated) problem. The second part
of the section is devoted to error estimates for time-incremental approximations, which are often used in numerical analysis of evolution
problems. The third part of Section 3 concerns the estimate of deviations from the exact solution to the parabolic thin obstacle problem. Finally, in Section~4, we consider several examples that demonstrate how the  estimates work in practice.


\section{Estimates of deviations from the exact solution
to the parabolic obstacle problem}
\subsection{Problem setting}
We consider the classical parabolic obstacle problem, which elliptic part
is presented by the Laplace operator. For simplicity, we restrict our consideration to the case of  time-independent obstacles.

Let 
$\Omega$ be an open, connected, and bounded domain in $\mathbb{R}^n$ with Lipschitz continuous boundary $\partial\Omega$, $Q_T=\Omega \times ]0,T[$.
We consider an obstacle function $\phi$ satisfying
$$
\phi \in H^2(\Omega)  \qquad \text{and}\qquad \phi \leqslant 0 \quad \text{a.e. on}\ \, \partial\Omega.
$$
 The class of admissible functions is defined as follows:
$$
\mathbb{K}=\mathbb{K}(\phi):=\{ w\in L^2((0,T),H^1_0(\Omega)):\  w_t\in L^2((0,T),H^{-1}(\Omega)), \ w \geqslant \phi(x) \ \; \text{a.e. on}\ \, Q_T\}.
$$
By a standard interpolation argument, the above assumptions imply $w\in C^0((0,T),L^2(\Omega))$.
Note that $K$ is non-empty due to the compatibility condition $\phi \leqslant 0$ on the lateral boundary $\partial\Omega \times (0,T)$ (which has to be understood in the sense of traces).
Henceforth, we assume that $f\in L^2(Q_T)$ and
$$
u_0\in H^1_0(\Omega) \quad \text{with}\quad u_0\geqslant \phi\quad \text{a.e. on}\ \Omega.
$$





We consider the following variational {\it Problem $\cP$.}\;
Find a function 
$u\in \mathbb{K}$
such that for almost all $t$ and $\forall w\in \mathbb{K}$ we have
\begin{eqnarray}
&&\int\limits_{Q_T}u_t (w-u)dxdt+\int\limits_{Q_T}\nabla u \cdot \nabla (w-u)dxdt \geqslant \int\limits_{Q_T}f(w-u)dxdt, 
\label{eq:eq1}\\
&& u(x,0)=u_0(x), \quad \forall x\in \Omega.
\label{eq:bc2}
\end{eqnarray}
Here and later on, $w_t$, resp. $\frac{\partial w}{\partial t}$ denotes the partial derivative with respect to time
and  $\nabla w$ denotes the spatial gradient vector.

It is known (see, e.g., \cite{LS1967}, \cite{B1972}, \cite{B1972a}, and \cite{DLi1976}) that under the above assumptions the minimization  problem (\ref{eq:eq1})-(\ref{eq:bc2}) is uniquely solvable.

By $Q^+_{T}(u):=\{(x,t)\in Q_T \mid u(x,t)>\phi\}$ we denote the subset
of $Q_T$, where the obstacle is not active. In this set,
\begin{equation}
\label{obs1}
f+\dvg p-u_t= 0,\qquad p=\nabla u.
\end{equation}
In the reminder (coincidence) set
$Q^\phi_{T}(u):=\{(x,t)\in Q_T \mid u(x,t)=\phi\}$
it holds
\begin{equation*}
f+\dvg p-u_t\leq 0.
\end{equation*}

\begin{remark} It is also well known (see, e.g., \cite{ASU2000}, \cite{CPS2004}, or \cite{A2018}) that the best possible regularity of a solution, $u$, to a parabolic obstacle problem is $u\in W^{2,1,\infty}_{loc}(Q_T)$, even when the source, boundary data, and obstacle function and domain boundary are $C^{\infty}$.
\end{remark}

\subsection{Estimates of the distance to the exact solution}
Let $v\in \mathbb{K}$ be a function viewed as an approximation
of the exact solution $u$, so that $e:=v-u$ is the error
and
$Q^+_{T}(v):=\{(x,t)\in Q_T \mid v(x,t)>\phi\}$ and
$Q^\phi_{T}(v):=\{(x,t)\in Q_T \mid v(x,t)=\phi\}$ denote the sets
associated with $v$.
Our goal is to deduce a computable majorant of $e$, which uses only known
information (i.e., the function $v$, sets $Q^+_{T}(v)$ and $Q^\phi_{T}(v)$,
$u_0$, $\Omega$, and other data of the Problem $\cP$). The error is measured
in terms of the combined error norm
\begin{equation}
|[e]|^2_{\alpha,Q_T}:=\|e(\cdot,T)\|^2_{\Omega}+\left(2-\frac{1}{\alpha}\right)\|\nabla e\|^2_{Q_T},\quad \alpha \geqslant \dfrac{1}{2}.
\label{eq:3.1}
\end{equation}

For this purpose we combine the methods earlier developed
for stationary problems with obstacles (see \cite{Re2000}, \cite{Re2007}, \cite{Re2008}, \cite{AR2018}, \cite{AR2020})
and for
parabolic equations (see \cite{Re2002}, \cite{MR2016}, \cite{LMR2019}).

In the above cited publications, the reader can also find numerical examples conforming the efficiency of these estimates to
problems with obstacles and finite element and IgA approximations of
evolutionary problems.

\begin{theorem}
\label{Th1}
For any $v\in {\mathbb K}$ and any vector valued function $\tau$ such that
 \begin{equation*}
 \tau  \in H_{\dive }(Q_T):=\left\{\tau  (x,t)\in L^2(Q_T, \mathbb{R}^n) \mid \dive\tau  \in L^2(\Omega) \ \text{for a.e.}\ t\in (0,T) \right\}
 \end{equation*}
  it holds
\begin{equation}
|[e]|^2_{\alpha,Q_T} \leqslant \|e(\cdot,0)\|^2_{\Omega} +\alpha \left( \|\tau -\nabla v\|_{Q_T}+C_F\|\mathcal{F}_f(v,\tau)\|_{Q_T}\right)^2,
\label{eq:mainest}
\end{equation}
where
$$
\mathcal{F}_f(v, \tau ):=\left\{\begin{aligned}
&\mathcal{R}_f(v, \tau ), \  &&\text{if}\ \, (x,t) \in Q^+_{T}(v),\\
&\{\mathcal{R}_f(v, \tau)\}_{\oplus}, \  &&\text{if}\ \, (x,t) \in Q^\phi_T(v),
\end{aligned}
\right.
$$
and $
\mathcal{R}_f(v,\tau ):=f+\dive\tau  - v_t$. The right hand side
of (\ref{eq:mainest}) vanishes if an only if $v=u$ and $\tau=\nabla u$.
\end{theorem}
\begin{proof}
From (\ref{eq:eq1}) it follows that
\begin{align*}
\int\limits_{Q_T} (u_t-v_t) (w-u)&dxdt+\int\limits_{Q_T}\nabla (u-v) \cdot \nabla (w-u)dxdt  \\
&\geqslant\int\limits_{Q_T} f(w-u)dxdt -\int\limits_{Q_T} v_t (w-u)dxdt -\int\limits_{Q_T} \nabla v \cdot \nabla (w-u)dxdt,
\end{align*}
for any $w\in \mathbb{K}$. In particular, for $w=v$ we have
\begin{equation} \label{eq:3.3}
\begin{aligned}
\frac{1}{2} \int\limits_{Q_T} \frac{\partial (u-v)^2}{\partial t}dxdt&+\int\limits_{Q_T}|\nabla (u-v)|^2dxdt
\leqslant \int\limits_{Q_T}f(u-v)dxdt \\&-\int\limits_{Q_T} v_t  (u-v)
-\int\limits_{Q_T}\nabla v \cdot \nabla (u-v)dxdt.
\end{aligned}
\end{equation}

Since $\tau  \in H_{\dive }(Q_T)$ and 
$u-v\in H^1_0(\Omega)$, the identity
\begin{equation} 
\label{eq:3i}
\int\limits_{\Omega} \tau  \cdot \nabla (u-v)dx=-\int\limits_{\Omega} (u-v) \dive \tau  dx
\end{equation}
holds for almost all $t\in (0,T)$.
Notice that
$$
\int\limits_{Q_T} \frac{\partial (u-v)^2}{\partial t}dxdt=\int\limits_{\Omega} (u-v)^2 dx \bigg|^T_0.
$$
Hence using definitions of $\mathcal{R}_f(v,\tau )$ and $\mathcal{F}_f(v,\tau )$,
we write (\ref{eq:3.3})  in the form
\begin{equation} \label{eq:3.4}
\begin{aligned}
\frac{1}{2} \|(u-v)(\cdot,T)\|^2_{\Omega} &-\frac{1}{2} \|(u-v)(\cdot,0)\|^2_{\Omega} +\|\nabla (u-v)\|^2_{Q_T}\\
&\leqslant \int\limits_{Q_T} \left( \mathcal{R}_f(v,\tau )\right)(u-v)dxdt +
\int\limits_{Q_T}\left( \tau  - \nabla v\right)\nabla (u-v)dxdt\\
&\leqslant \int\limits_{Q_T} \mathcal{F}_f(v,\tau )(u-v)dxdt + \int\limits_{Q_T}\left( \tau  - \nabla v \right)\nabla (u-v)dxdt.
\end{aligned}
\end{equation}

We set $e:=u-v$. Estimating the first term on the right-hand side of (\ref{eq:3.4})  by the Friedrich's type inequality and the second term there by the H{\"o}lder inequality,  we arrive at
\begin{equation}
 \label{eq:3.5}
\begin{aligned}
\frac{1}{2}\|e(\cdot, T)\|^2_{\Omega}+\|\nabla e\|^2_{Q_T} &\leqslant \frac{1}{2}\|e(\cdot, 0)\|^2_{\Omega}
+\|\tau  - \nabla v\|_{Q_T} \|\nabla e\|_{Q_T}\\
&+C_F \|\mathcal{F}_f(v,\tau ) \|_{Q_T} \|\nabla {e}\|_{Q_T}.
\end{aligned}
\end{equation}

Since
\begin{equation}
\label{eq:3.5a}
\begin{aligned}
\bigg( \|\tau  - \nabla v\|_{Q_T}&+C_F \|\mathcal{R}_f(v,\tau )
 \|_{Q_T} \bigg) \|\nabla {e}\|_{Q_T}\\
& \leqslant 
\frac{\alpha}{2} \bigg( \|\tau  - \nabla v\|_{Q_T}+C_F \|\mathcal{F}_f(v,\tau ) \|_{Q_T} \bigg)^2+\frac{1}{2\alpha}
\| \nabla {e}\|^2_{Q_T},
\end{aligned}
\end{equation}
the inequality (\ref{eq:3.5}) yields the estimate
\begin{equation} \label{eq:3.6}
\begin{aligned}
\frac{1}{2}\|e(\cdot, T)\|^2_{\Omega}+ \left(1-\frac{1}{2\alpha}\right) \| \nabla {e}\|^2_{Q_T} &\leqslant 
\frac{1}{2}\|e(\cdot, 0)\|^2_{\Omega} \\
&+\frac{\alpha}{2} \bigg(\|\nabla v-\tau \|_{Q_T} +C_F\|\mathcal{F}(v,\tau )\|_{Q_T}\bigg)^2.
\end{aligned}
\end{equation}
Now (\ref{eq:mainest}) follows from (\ref{eq:3.6}) provided that the condition $\|e(\cdot,0)\|^2_{\Omega}=0$ is fulfilled (which is obviously the case for $v\in \mathbb{K}$).

Assume that the right hand side of (\ref{eq:mainest}) vanishes. Then
$\tau=\nabla v$, $v(x,0)=u_0(x)$ and
\begin{align}
\mathcal{R}_f(v,\tau )&=0\quad{\rm in}\;Q^+_T(v), \notag\\
\label{eq:Rf2}
\mathcal{R}_f(v,\tau )&\leq 0\quad{\rm in}\;Q^\phi_T(v).
\end{align}
We use these relations to estimate the integral
\begin{multline} 
\label{eq:unique}
\int\limits_{Q_T}\left(
v_t (w-v)dxdt+\nabla v \cdot \nabla (w-v)-f(w-v)
\right)dxdt\\=\int\limits_{Q^+_T(v)}\left(
v_t (w-v)dxdt+\tau \cdot \nabla (w-v)-f(w-v)
\right)dxdt\\
+
\int\limits_{Q^\phi_T(v)}\left(
v_t (w-v)dxdt+\tau \cdot \nabla (w-v)-f(w-v)
\right)dxdt\\
=\int\limits_{Q^\phi_T(v)}\left(
(\dvg\tau +f-v_t)(v-w)
\right)
dxdt.
\end{multline}
In view of (\ref{eq:Rf2}), both  terms of the integrand in the right-hand side of (\ref{eq:unique})
are nonpositive for any $w\in {\mathbb K}$. Hence we see that $v$ satisfies (\ref{eq:eq1}). Since the solution is unique, we conclude that $v$ coincides with $u$ and $\tau$ coincides with the exact flux $p=\nabla u$.
\end{proof}

It is worth adding some comments to Theorem \ref{Th1}.
\begin{remark}
 The left hand
side of (\ref{eq:mainest}) is a natural measure of the distance between
$v$ and $u$, which particular form depend on the parameter $\alpha$. The left hand side is directly computable.
Since $\alpha$ appear as a multiplier in the right hand side, we should not select it too large. For $\alpha=1$, we obtain a useful estimate
\begin{equation}
\label{eq:spec1}
\|e(\cdot,T)\|^2_{\Omega}+\|\nabla e\|^2_{Q_T}\leq
\|e(\cdot,0)\|^2_{\Omega} + \left(\|\tau -\nabla v\|_{Q_T}+C_F\|\mathcal{F}_f(v,\tau)\|_{Q_T}\right)^2.
\end{equation}
Another simple estimate corresponds to the limit case $\alpha=\frac12$:
\begin{equation}
\label{eq:spec2}
\|e(\cdot,T)\|^2_{\Omega}\,\leq
\|e(\cdot,0)\|^2_{\Omega} + \frac{1}{2}\left(\|\tau -\nabla v\|_{Q_T}+C_F\|\mathcal{F}_f(v,\tau)\|_{Q_T}\right)^2.
\end{equation}
We outline that the estimates (\ref{eq:mainest}), (\ref{eq:spec1}),
and (\ref{eq:spec2}) are valid for {\em any} function $v\in {\mathbb K}$
regardless of the method by which it is constructed. This is the principal difference of functional type a posteriori estimates from
other a posteriori estimates, which usually impose special conditions
on approximations (e.g., Galerkin orthogonality). In the next section, we
use this universality feature to deduce simple bounds of errors caused
by simplifications of the source term and initial condition.
\end{remark}

\begin{remark}
A particular form of (\ref{eq:mainest}) can be viewed as a generalisation
of the well known hypercircle estimate to the case of the parabolic
obstacle problem.
Define the set
\be
Q_{f,\phi}:=\{\tau\in L^2(\Omega,\Rn): \mathcal{R}_f(v, \tau )=0\;\text{in}\; Q^+_{T}(v), \;\mathcal{R}_f(v, \tau )\leq 0\;\text{in}\; Q^\phi_{T}(v)\}
\ee
We have
\begin{equation*}
|[e]|^2_{\alpha,Q_T} \leqslant \|e(\cdot,0)\|^2_{\Omega} +\alpha^2 \|\tau -\nabla v\|^2_{Q_T}
\end{equation*}
for any pair of functions $v$ and $\tau$ such that 
$(v,\tau)\in Q_{f,\phi}$. Notice that unlike the estimates known for linear problems (which contain only the restriction $\dvg \tau+f=0$
for the dual variable), this estimate impose the condition on 
both functions $v$ and $\tau$. This effect is generated by the nonlinearity associated with existence of a coincidence set and free
boundary.
\end{remark}
\section{Special cases}
\subsection{Errors generated by simplification of the model}

Simplification (coarsening, defeaturing) of a mathematical model may be very useful
if we can eliminate insignificant details without essential loss of the
accuracy. Simplification methods for elliptic type problems are
well studied (the reader can find a complete theory in \cite{RS2020}).
Here we briefly discuss these questions in the context of Problem $\cP$.

Consider Problem $\wt\cP$, which uses a function $\widetilde{f}(x,t)$ instead of $f(x,t)$ and $\widetilde{u}_0(x)$ instead of $u_0(x)$. Let 
$\widetilde{u}(x,t)$ and $\widetilde{p}=\nabla \widetilde{u}$
be the corresponding exact solution and exact flux. Substituting these functions in (\ref{eq:3.1}), we obtain a simple estimate
\begin{equation*} 
|[u-\widetilde{u}]|^2_{\alpha, Q_T} \leqslant \|u_0-\widetilde{u}_0\|^2_{\Omega}+\alpha C_F\|\mathcal{F}(\widetilde{u}, \widetilde{p})\|^2_{Q_T}.
\end{equation*}
If $(x,t) \in Q^+_T(\widetilde{u})$ then
$$
\mathcal{R}_f(\widetilde{u}, \widetilde{p})=f+  \dive\,\widetilde{p}-\partial_t \widetilde{u}=f-\widetilde{f}.
$$
If $(x,t)\in Q^\phi_T(\widetilde{u})$, then
$$
\mathcal{R}_{\wt f}(\widetilde{u}, \widetilde{p})=\widetilde{f}+\dive\,\widetilde{p}-\partial_t\widetilde{u} \leqslant 0
$$
and we find that
$$
\max \{0, f+\dive\,\widetilde{p}-\partial_t \widetilde{u}\} 
=\max \{0,f-\widetilde{f}+\mathcal{R}_{\wt f}(\widetilde{u}, \widetilde{p})\} \leqslant \max \{0, f-\widetilde{f}\}=\{f-\widetilde{f}\}_{\oplus}.
$$
Hence (\ref{eq:3.3}) yields the estimate
\begin{equation} \label{eq:simple}
|[u-\widetilde{u}]|^2_{\alpha, Q_T} \leqslant \|u_0-\widetilde{u}_0\|^2_{\Omega}+\alpha C_F \|g(x,t)\|^2_{Q_T},
\end{equation}
where
$$
g(x,t)=\left\{
\begin{aligned}
&f-\widetilde{f}, && \text{if}\ (x,t)\in Q^+_T(\widetilde{u}),\\
&\{f-\widetilde{f}\}_{\oplus}, && \text{if}\  (x,t)\in Q^\phi_T(\widetilde{u}).
\end{aligned}
\right.
$$
Since $\|g(x,t)\|_{Q_T} \leqslant \|f-\widetilde{f}\|_{Q_T}$, we have
a simplified estimate
\begin{equation} \label{eq:sim}
|[u-\widetilde{u}]|^2_{\alpha,Q_T} \leqslant \|u_0-\widetilde{u}_0\|^2_{\Omega}+\alpha C_F \|f-\widetilde{f}\|^2_{Q_T}.
\end{equation}
In general, the estimate (\ref{eq:sim}) is coarser than (\ref{eq:simple}), but it does not require knowledge on the coincidence set $Q^\phi_T(\widetilde{u})$.
\subsection{Errors of time--incremental approximations}

Now we consider a special class of approximations, which are typically used in time-incremental methods for various evolutionary problems.
Let the interval $(0,T)$ be split into a collection of subintervalls
\begin{equation*} 
I_k:=(t_k,t_{k+1}), \quad t_{k+1}-t_k=\Delta_k>0, \quad t_0=0, \quad t_N=T,
\end{equation*}
and the approximation $v(x,t)$ has the form
\begin{equation} \label{eq:increment}
v(x,t)=v_k(x)+\frac{v_{k+1}(x)-v_k(x)}{\Delta_k}(t-t_k)
\end{equation}
for $(x,t)\in Q_k:=\Omega \times I_k$. Here $v_k(x)\in H^1_0(\Omega)$, $v_k(x) \geqslant \phi$ are the approximations computed by a time-incremental 
numerical method. 
Notice that so defined function $v$ belongs to ${\mathbb K}$. Indeed 
\be
v(x,t)-\phi(x)=v_k(x)\frac{t_{k+1}-t}{\Delta_k}+
v_{k+1}(x)\frac{t-t_k}{\Delta_k}-\phi(x)\geq 0.
\ee
This property also holds in more complicated cases, where $\phi$ depends
on  on time (e.g., if $\phi$ is linear function of $t$). However, for simplicity
in this paper we consider only the case $\phi=\phi(x)$.

We define the sets
$$
\Omega^\phi_k(v_k):=\{x\in \Omega : v_k(x)=\phi\}\qquad\text{and}\quad \Omega^+_k(v_k):=\{x\in \Omega : v_k(x)>\phi\}.
$$
 Notice that for the function $v$ defined by (\ref{eq:increment}), the
set $Q^\phi_k:=Q_k\cap Q^\phi(v)$ is defined as follows:
\begin{equation*} 
Q^\phi_k(v):=Q^\phi_T(v)\cap Q_k=\Omega^{\phi}_{k+1/2}\times I_k.
\end{equation*}
where $\Omega^{\phi}_{k+1/2}:=\Omega_k^\phi \cap \Omega^ \phi_{k+1}$.
 For $(x,t)\in Q_k$ we have
\begin{equation*} 
v_t =\frac{1}{\Delta_k}(v_{k+1}-v_k), \qquad \nabla v=\nabla v_k+\frac{\nabla(v_{k+1}-v_k)}{\Delta_k}(t-t_k).
\end{equation*}

Below we deduce two different a posteriori estimates for semi-discrete approximation. The simplest estimate is valid for the case, where $\Delta_k$ is so small, that we can neglect changes of the source term $f(x,t)$ and replace  it by the averaged function
\begin{equation} \label{eq:4.2}
\widetilde{f}(x,t)=\mean{f}_{I_k}(x):=\frac{1}{\Delta_k}\int\limits_{I_k}f(x,t)dt.
\end{equation}
In this situation,  it is natural to select for the flux the simplest approximation also:
\begin{equation} \label{eq:4.3}
\tau (x,t)=\tau_k(x) \quad \text{for}\quad \;(x,t)\in Q_k.
\end{equation}
A more advanced version uses an affine approximation of $f(x,t)$:
\begin{equation} \label{eq:4.4}
\widetilde{f}(x,t)=f_k(x)+\frac{f_{k+1}(x)-f_k(x)}{\Delta_k}(t-t_k)+\zeta_k(x) \quad \text{for}\quad (x,t)\in Q_k,
\end{equation}
where $\zeta_k(x)=\mean{f}_{I_k}(x)-\frac12(f_k(x)+f_{k+1}(x))$. It is selected such that
$$
\int\limits_{I_k}\widetilde f dt=
\int\limits_{I_k}f(x,t) dt=\Delta_k\mean{f}_{I_k}(x).
$$

A similar time--incremental
form can be used to approximate the flux $\tau(x,t)$ in $Q_k$. Let $\tau_k(x)$, $k=0,1,2,...$ be approximations related to $t_k$.
We set
\begin{equation} \label{eq:4.5}
\tau (x,t)=\tau_k(x)+\frac{\tau_{k+1}(x)-\tau_k(x)}{\Delta_k} (t-t_k).
\end{equation}

Let $\widetilde{u}$ be the exact solution of the problem $\wt\cP$ (where $f$ is replaced by $\widetilde{f})$.  We have
\begin{equation*}
|[u-v]|_{\alpha,Q_T} \leqslant |[u-\widetilde{u}]|_{\alpha,Q_T}+|[\widetilde{u}-v]|_{\alpha,Q_T}.
\end{equation*}

Here  the first term in the right-hand side is estimated by (\ref{eq:sim}),
which contains only the second term (because in our case $\wt u_0=u_0$). Hence we need to estimate the last term only.
For this purpose, we use (\ref{eq:3.4}) and obtain
\begin{equation} \label{eq:4.6}
\begin{aligned}
2\|\nabla (v-\widetilde{u})\|^2_{Q_T}&+\|(v-\widetilde{u})(\cdot,T)\|^2_{\Omega} \leqslant \|v_0-u_0\|^2_{\Omega}\\
&+2\sum\limits_{k=0}^{N-1} \int\limits_{Q_k}\left[(\tau -\nabla v)\nabla (\widetilde{u}-v)+\mathcal{F}_{\widetilde{v}}(v,\tau )(\widetilde{u}-v)\right]dxdt.
\end{aligned}
\end{equation}
Notice that
$$
\int\limits_{Q_k}\mathcal{F}_{\widetilde{f}}(v,\tau )(\widetilde{u}-v)dxdt=\int\limits_{I_k}\int\limits_{\Omega}\mathcal{F}^k_{\widetilde{f}}(v,\tau )(\widetilde{u}-v)dxdt,
$$
where
$$
\mathcal{F}^k_{\widetilde{f}}(v,\tau ):=\left\{
\begin{aligned}
&\mathcal{R}_{\widetilde{f}}(v,\tau ) &&\text{if}\quad x\in \Omega \setminus \Omega^\phi_{k+1/2},\\
\{&\mathcal{R}_{\widetilde{f}}(v,\tau )\}_{\oplus} &&\text{if}\quad x\in  \Omega^\phi_{k+1/2}.
\end{aligned}
\right.
$$

Consider first the simplest estimate that follows from (\ref{eq:4.6}) with $\widetilde{f}$ and $\tau$ selected in accordance with (\ref{eq:4.2}) and (\ref{eq:4.3}), respectively. In this case,
$$
\mathcal{F}^k_{\widetilde{f}}(v,\tau )=\left\{
\begin{aligned}
&\mathcal{R}^k_{\widetilde{f}}(v_k,v_{k+1}, \tau_k) &&\text{if}\quad x\in \Omega \setminus \Omega^\phi_{k+1/2},\\
\{&\mathcal{R}^k_{\widetilde{f}}(v_k,v_{k+1}, \tau_k)\}_{\oplus} &&\text{if}\quad x\in \Omega^\phi_{k+1/2}, 
\end{aligned}
\right.
$$
where
\begin{equation*} 
\mathcal{R}^k_{\widetilde{f}}(v_k,v_{k+1}, 
\tau_k):=\widetilde{f}(x)+\dive\,\tau_k(x)-\frac{v_{k+1}-v_k}{\Delta_k}
\end{equation*}
depends on $x$ only. Therefore,
\begin{equation} \label{eq:4.7}
\begin{aligned}
\int\limits_{Q_k}\mathcal{F}^k_{\widetilde{f}}(v,\tau )(\widetilde{u}-v)dxdt & \leqslant C_F\|\mathcal{F}^k_{\widetilde{f}}(v,\tau_k )\|_{\Omega}
\int\limits_{I_k}\|\nabla (\widetilde{u}-v)\|_{\Omega}dt\\
&\leqslant C_F  \Delta^{1/2}_k \|\mathcal{F}^k_{\widetilde{f}}(v,\tau_k )\|_{\Omega} \|\nabla (\widetilde{u}-v)\|_{Q_k}.
\end{aligned}
\end{equation}
Next,
\begin{equation*} 
\int\limits_{Q_k}(\tau -\nabla v)(\widetilde{u}-v)dxdt \leqslant \|\tau -\nabla v\|_{Q_k}\|\nabla(\widetilde{u}-v)\|_{Q_k},
\end{equation*}
where
\begin{equation*} 
\nabla v-\tau=\frac{\nabla (v_{k+1}-v_k)}{\Delta_k}(t-t_k)+\nabla v_k -\tau_k.
\end{equation*}
Let
$$
D_1^k(v_k,v_{k+1}):=\frac{1}{12}\|\nabla (v_{k+1}-v_k)\|^2_{\Omega}
$$
and
$$
D_2^k(v_k,v_{k+1},\tau_k):=\|\frac{1}{2}\nabla (v_{k+1}+v_k)-\tau_k\|^2_{\Omega}.
$$
Then,
\begin{equation} \label{eq:4.8}
\int\limits_{Q_k}|\tau -\nabla v|^2dxdt=\int\limits_{\Omega}\int\limits_{I_k} |\tau -\nabla v|^2dtdx=\Delta_k\Bigl(D_1(v_k,v_{k+1})+D_2(v_k,v_{k+1}, \tau_k)\Bigr).
\end{equation}

By (\ref{eq:4.6}), (\ref{eq:4.7}), and (\ref{eq:4.8}), we obtain
\begin{align*}
2\|\nabla (\widetilde{u}-v)\|^2_{Q_T} &+ \|(\widetilde{u}-v)(\cdot, T)\|^2_{\Omega} \leqslant \|u_0-v_0\|^2_{\Omega}\\&+
2\sum\limits_{k=0}^{N-1}\Delta^{1/2}_k \left[
(D_1^k+D_2^k)^{1/2}+C_F\|\mathcal{F}^k(v, \tau )\|_{\Omega}\right]\|\nabla(\widetilde{u}-v)\|_{Q_k}\\
&
\leqslant \|u_0-v_0\|^2_{\Omega}\\
&+2\|\nabla (\widetilde{u}-v)\|_{Q_T}\Bigl(\sum\limits_{k=0}^{N-1}\Delta_k \left[(D_1^k+D_2^k)^{1/2}+{C_F}\|\mathcal{F}^k(v,\tau )\|_{\Omega}\right]^2\Bigr)^{1/2}.
\end{align*}

After using Young's inequality, we arrive at the estimate
\begin{equation} \label{eq:4.9}
|[\widetilde{u}-v]|^2_{\alpha,Q_T} \leqslant \|u_0-v_0\|^2_{\Omega}+\alpha \left(\sum\limits_{k=0}^{n-1}
\Delta_k \left[(D_1^k+D_2^k)^{1/2}+C_F \|\mathcal{F}^k(v, \tau )\|_{\Omega}\right]^2\right)^{\frac{1}{2}}.
\end{equation}

Now (\ref{eq:3.4}) and (\ref{eq:4.9}) imply the desired error majorant
\begin{equation} \label{eq:4.10}
\begin{aligned}
|[u-v]|^2_{\alpha,Q_T} &\leqslant \|u_0-v_0\|^2_{\Omega}+\alpha C_F\|f-\widetilde{f}\|^2_{Q_T}\\&+\alpha \left(\sum\limits_{k=0}^{n-1}
\Delta_k \left[(D_1^k+D_2^k)^{1/2}+C_F \|\mathcal{F}^k(v, \tau )\|_{\Omega}\right]^2\right)^{\frac{1}{2}}.
\end{aligned}
\end{equation}

The first two terms in the right hand side of (\ref{eq:4.10}) reflect errors generated by simplification of the initial data and source term. The last term reflects the errors caused by semi-discrete approximations. If the approximation is sufficiently regular and has no jumps, then the term $D_1(v_k,v_{k+1})$ is of the order $\Delta^2_k$, i.e. it is a minor term. The main term is $D_2(v_k, v_{k+1},\tau_k)$. It penalises inaccuracy in the relation $p^*=\nabla u$, which must hold for the exact solution and its flux. This term is small if the mean gradient $\nabla \left(\frac{v_{k+1}+v_k}{2}\right)$ is close to the flux approximation $\boldsymbol{\sigma}_k$ in $Q_k$. The functions $\boldsymbol{\sigma}_k$ can be viewed as images of the quantity $\frac{1}{\Delta_k}\int\limits_{I_k}p(x,t)dt$ associated with the true flux $p=\nabla u$. They can be extracted from the numerical solution $v$. 

For example, let $\Omega$ be a poligonal domain discretized by a simplicial mesh $\mathcal{F}_h$ and $v^h_k$, $k=0,1,\dots,N$ denote the respective numerical solutions computed using the finite element method 
for each step of the time-incremental sequence. Well known gradient averaging methods generate ``averaged'' fluxes
\begin{equation*} 
\boldsymbol{\sigma}^h_k(x)=G_k\nabla v^h_k \in H_{\dive}(\Omega),
\end{equation*}
where $G_h$ is an averaging operator. Such an operator
can be based on a simple  patch-averaging procedure
 or use more complicated procedures of global averaging (see, for example,
 \cite{CaBa2002} and \cite{BaCa2004}). 

Then, we can set
\begin{equation*} 
\tau_k=\frac{\boldsymbol{\sigma}^h_k+\boldsymbol{\sigma}^h_{k+1}}{2}.
\end{equation*}
Moreover, we are not limited to such a choice of $\tau_k$, $k=1,2,\dots, N$, which may be regarded as an initial guess only. The majorant (\ref{eq:4.10}) allows us to modify them this function
in order to minimise the right hand side of (\ref{eq:4.10}).

Now, we consider a more advanced error bound, which follows from (\ref{eq:4.4}), (\ref{eq:4.5}), and (\ref{eq:4.6}). In this case
$$
\mathcal{R}^k_{\widetilde{f}}(v,\tau )=
\mathcal{R}^k_{f_k}(v_k,v_{k+1},\tau_k)+\frac{\mathcal{R}^k_{f_{k+1}}-\mathcal{R}^k_{f_k}}{\Delta_k}(t-t_k),
$$
where
\begin{align*}
\mathcal{R}^k_{f_k}(v_k,v_{k+1}, \tau_k)&=f_k+\dive\,\boldsymbol{\sigma}_k-\frac{v_{k+1}-v_k}{\Delta_k},\\
\mathcal{R}^k_{f_{k+1}}(v_k,v_{k+1}, \boldsymbol{\sigma}_{k+1})&=f_{k+1}+\dive\,\boldsymbol{\sigma}_{k+1}-\frac{v_{k+1}-v_k}{\Delta_k}
\end{align*}
and define the sets
\begin{align*}
\omega^k_1&:=\left\{x\in \Omega^{\phi}_{k+1/2} : \mathcal{R}^k_{f_k}(v_k, v_{k+1},\tau_k) \leqslant 0 \right\},\\
\omega^k_2&:=\left\{x\in \Omega^{\phi}_{k+1/2} : \mathcal{R}^k_{f_{k+1}}(v_k, v_{k+1},\boldsymbol{\sigma}_{k+1}) \leqslant 0 \right\},\\
\omega^k&:=\omega^k_1 \cap \omega^k_2.
\end{align*}

The \ functions \ $\mathcal{R}^k_{f_k}(v_k, v_{k+1},{\tau}_k)$ \ and \ $\mathcal{R}^k_{f_{k+1}}(v_k, v_{k+1},{\tau}_{k+1})$ \ have \ clear \ meanings. \ They \ present residuals of the differential equation (in the incremental form, where time derivative is replaced by the finite difference) associated with the boundary points $t_k$ and $t_{k+1}$ of the interval$I_k$.

It is not difficult to see that 

\begin{equation} \label{eq:4.11}
\begin{aligned}
\int\limits_{I_k}\int\limits_{\Omega}\mathcal{F}^k_{\widetilde{f}}(v,\tau )(\widetilde{u}-v)dxdt &\leqslant \int\limits_{I_k}\int\limits_{\Omega \setminus \omega^k} \mathcal{F}^k_{\widetilde{f}}(v,\tau )(\widetilde{u}-v)dxdt \\
&\leqslant C_F \|\mathcal{F}^k_{\widetilde{f}}(v,\tau )\|_{I_k \times (\Omega \setminus \omega^k)}\|\nabla (\widetilde{u}-v)\|_{Q_k}.
\end{aligned}
\end{equation}
Here

\begin{equation}\label{eq:4.12}
\begin{aligned}
\|\mathcal{F}^k_{\widetilde{f}}(v,\tau )\|_{I_k \times (\Omega \setminus \omega^k)}^2&=\int\limits_{\Omega \setminus \omega^k}\int\limits_{I_k}|\mathcal{F}^k_{\widetilde{f}}(v,\tau )|^2dtdx \\
&=
\int\limits_{\Omega \setminus \omega^k}\int\limits_{t_k}^{t_{k+1}}\left(
\mathcal{R}^k_{f_k}+\frac{\mathcal{R}^k_{f_{k+1}}-\mathcal{R}^k_{f_k}}{\Delta_k}(t-t_k)\right)^2dtdx\\
&=\frac{\Delta_k}{4}\int\limits_{\Omega \setminus \omega^k}
\left[\left(\mathcal{R}^k_{f_{k+1}}-\mathcal{R}^k_{f_k}\right)^2 +\frac{1}{3}\left(\mathcal{R}^k_{f_{k+1}}+\mathcal{R}^k_{f_k}\right)^2
\right]dx\\
&=\frac{\Delta_k}{4}\left[
\|\mathcal{R}^k_{f_{k+1}}-\mathcal{R}^k_{f_k}\|_{\Omega \setminus \omega^k}^2+\frac{1}{3}\|\mathcal{R}^k_{f_{k+1}}+\mathcal{R}^k_{f_k}\|^2_{\Omega \setminus \omega^k}
\right].
\end{aligned}
\end{equation}

Consider another term. We have
\begin{equation*}
\int\limits_{Q_k}|\tau -\nabla v|^2dxdt=\frac{\Delta_k}{4} \left[
\|D^{k+1}-D^k\|^2_{\Omega}+\frac{1}{3}\|D^{k+1}+D^k\|^2_{\Omega}
\right],
\end{equation*}
where $D^k:=\tau _k-\nabla v_k$. By (\ref{eq:4.6}), (\ref{eq:4.11}), and (\ref{eq:4.12}), we obtain the following estimate:

\begin{equation}\label{eq:4.13}
\begin{aligned}
|[u-v]|^2_{\alpha, Q_T} &\leqslant \|u_0-v_0\|^2_{\Omega}+\alpha C_F \|\widetilde{f}-f\|^2_{Q_T}\\
&+\alpha \bigg(\sum\limits_{k=0}^{N-1} \frac{\Delta_k}{4} \bigg[ \|D^{k+1}-D^k\|^2_{\Omega}+\frac{1}{3}\|D^{k+1}+D^k\|^2_{\Omega} \\
 &+\|\mathcal{R}^k_{f_{k+1}}-\mathcal{R}^k_{f_k}\|_{\Omega \setminus \omega^k}^2
 +\frac{1}{3} \|\mathcal{R}^k_{f_{k+1}}+\mathcal{R}^k_{f_k}\|_{\Omega \setminus \omega^k}^2
\bigg]^{1/2}\bigg).
\end{aligned}
\end{equation}

The sum in the right hand side of (\ref{eq:4.13}) consists of the quantities that depend on the functions $v_k(x)$ and $\boldsymbol{\sigma}_k(x)$ that form semi-discrete approximations of the solution $u(x,t)$ and its flux $\nabla u(x,t)$. If the functions $D_k$, $\mathcal{R}^k_{f_k}$,  $\mathcal{R}^k_{f_{k+1}}$ are small (i.e., residuals of the relations generating the incremental posing of the problem), then (\ref{eq:4.13}) confirms accuracy of the computed solution.
\begin{remark}
A posteriori estimates for incremental approximations
of an evolutionary problem with obstacles has been studied
\cite{NSV2000} within the framework of the residual method, which operates with Galerkin approximations and uses special interpolation
operators related to the type of approximations selected. We use another conception, where the estimates are derived
on the functional level and, therefore, are independent of the
numerical method by which the approximation $v$ has been constructed.
The estimates do not contain local constants and use only the global constant $C_F$ associated with the domain $\Omega$.
At the same time, they contain a vector valued function $\tau$, which selection changes the majorant so that proper selection of this
function is important to have a sharp error estimate. Depending on the
method used it could be a difficult task (e.g., if we use standard
low order finite elements which usually produce rather coarse
approximations of fluxes) or a relatively simple one (e.g., if our
method generates mixed approximations and the corresponding fluxes
can be used  without post--processing). These issues require
further investigation. In simple examples presented below, we show
that realistic error bounds follow from the majorant even for very
simple reconstructions of the flux $\tau$. 
\end{remark}

\subsection{The parabolic Signorini problem}

A specific version of the problem $\mathcal{P}$ arises if an obstacle function $\phi$ is given on the part of the lateral surface of $Q_T$ instead of inside $Q_T$. 

Throughout of this subsection  we assume  that $f\in L^{\infty}(Q_T)$, $u_0\in W^2_{\infty}(\Omega)$, $\mathcal{M}$ is a relatively open subset of $\partial\Omega$ (in its relative topology), $\mathcal{S}=\partial\Omega \setminus \mathcal{M}$, and the obstacle function $\phi \in H^2(\mathcal{M})$ satisfies the compatibility conditions: 
$$
 \phi \leqslant u_0\ \text{on}\ \mathcal{M}, \quad  \phi \leqslant 0\ \text{on}\ \partial\mathcal{M}.
$$
Such a function $\phi$ is called the \textit{thin obstacle}.

In this case, the the problem (\ref{eq:eq1})--(\ref{eq:bc2}) is considered for   almost all $t\in (0,T)$ and all functions $w$ from the set 
\begin{align*}
\mathbb{K}_{\mathbb{S}}(\phi):=\{w\in L^2((0,T), H^1(\Omega)): w_t\in &L^2((0,T), H^{-1}(\Omega)), \\ &w\geqslant \phi\quad \text{a.e. on}\ \, \mathcal{M}_T, \quad w=0\quad \text{a.e. on}\ \, \mathcal{S}_T\},
\end{align*}
where $\mathcal{M}_T:=\mathcal{M}\times ]0,T[$, and $\mathcal{S}_T:=\mathcal{S}\times ]0,T[$.  This problem is known as the \textit{parabolic thin obstacle problem} or \textit{parabolic Signorini problem}.

Under the above assumptions on the problem data,  the existence of the unique solution has been established in \cite{LS1967}.  The exact solution $u$ satisfies the equations (\ref{obs1}) in $Q_T$ and the so-called \textit{Signorini boundary conditions}
$$
u \geqslant \phi,\quad \frac{\partial u}{\partial \mathbf{n}} \geqslant 0, \quad (u-\phi)\frac{\partial u}{\partial \mathbf{n}}=0 \quad \text{on}\ \, \mathcal{M}_T,
$$
where $\mathbf{n}$ denotes the unit outward normal to $\partial\Omega$. Moreover, according to \cite{ArU1988} and \cite{ArU1996}, the exact solution possesses H{\"o}lder continuous spatial gradient: $\nabla u \in \mathcal{H}^{\alpha, \alpha/2}_{\text{loc}}$ with the H{\"o}lder exponent $\alpha >0$ depending only on the dimension.

Furthermore, consider the set
\begin{align*}
H_{\dive}^{\mathbb{S}}(Q_T):=\{\tau (x,t)\in &L^2(Q_T, \mathbb{R}^n) \mid \dive\tau  \in L^2(\Omega), \\ &\tau\cdot \mathbf{n} \in L^2(\mathcal{M}),  \ \text{and}\ \tau \cdot \mathbf{n} \geqslant 0 \ \text{for a.e.}\ t\in (0,T) \}.
\end{align*}
Since for $v\in \mathbb{K}_{\mathbb{S}}$ we have $u-v =0$ a.e. on $\mathcal{S}_T$ only, the identity (\ref{eq:3i}) takes the form
$$
\int\limits_{\Omega} \tau \cdot \nabla (u-v)dx = -\int\limits_{\Omega}(u-v) \dive \tau dx	+\int\limits_{\mathcal{M}}(u-v)\,\tau\cdot \mathbf{n}d\mu
$$
for all $\tau \in H^{\mathbb{S}}_{\dive}(Q_T)$ and for almost all $t\in (0,T)$. 

Repeating all the arguments used in the derivation of (\ref{eq:3.4}) we conclude that 
\begin{equation}
\label{eq:33.0}
\begin{aligned}
\frac{1}{2}\|(u-v)(\cdot, T)\|^2_{\Omega}&-\frac{1}{2}\|(u-v)(\cdot,0)\|^2_{\Omega}+\|\nabla (u-v)\|^2_{Q_T}\\
&\leqslant \int\limits_{Q_T} \left(f+\dive\tau-v_t\right)(u-v)dxdt\\
&+\int\limits_{Q_T}(\tau-\nabla v)\nabla(u-v)dxdt
+\int\limits_{\mathcal{M}_T}(v-u)\tau\cdot \mathbf{n}d\mu dt.
\end{aligned}
\end{equation}

As in Section~2 we set $e:=u-v$ and estimate the first and second integrals on the right-hand side of (\ref{eq:33.0}) by the Friedrich's type inequality and the H{\"o}lder inequality, respectively. 
Notice that for  all $\tau \in H^{\mathbb{S}}_{\dive}(Q_T)$ and all $v\in \mathbb{K}_{\mathbb{S}}$ we can estimate the third integral on the right-hand side of (\ref{eq:33.0}) as follows:
$$
\begin{aligned}
\int\limits_{\mathcal{M}_T}(v-u)\tau\cdot \mathbf{n}d\mu dt&=
\int\limits_{\mathcal{M}_T}(v-\phi)\tau\cdot \mathbf{n}d\mu dt-
\int\limits_{\mathcal{M}_T}(u-\phi)\tau\cdot \mathbf{n}d\mu dt\\
&\leqslant \int\limits_{\mathcal{M}_T}(v-\phi)\,\tau\cdot \mathbf{n} \, d\mu dt.
\end{aligned}
$$
As a result we arrive at the estimate
\begin{equation}
\begin{aligned}
\frac{1}{2}\|e(\cdot,T)\|^2_{\Omega}+\|\nabla e\|^2_{Q_T}&\leqslant
\frac{1}{2}\|e(\cdot,0)\|^2_{\Omega}+\|\tau-\nabla v\|_{Q_T}\|\nabla e\|_{Q_T}\\
&+C_F\|f+\dive\tau-v_t\|_{Q_T}\|\nabla e\|_{Q_T}+\int\limits_{\mathcal{M}_T} (v-\phi)\, \tau \cdot \mathbf{n}\, d\mu dt.
\end{aligned}
\end{equation}

Taking into account (\ref{eq:3.5a}) we conclude that  the estimate
\begin{equation}
\label{eq: 33.1}
\begin{aligned}
\frac{1}{2}\|e(\cdot,T)\|^2_{\Omega}&+\left(1-\frac{1}{2\alpha}\right)\|\nabla e\|^2_{Q_T}\\
&\leqslant \frac{1}{2}\|e(\cdot,0)\|^2_{\Omega}+
\frac{\alpha}{2}\left(
\|\nabla v-\tau\|_{Q_T}+C_F\|\mathcal{F}(v,\tau)\|_{Q_T}\right)^2\\
& +\int\limits_{\mathcal{M}_T} (v-\phi) \, \tau \cdot \mathbf{n}\, d\mu dt
\end{aligned}
\end{equation}
holds true for any $\alpha \geqslant \frac{1}{2}$, all $v\in \mathbb{K}_{\mathbb{S}}$,  and all $\tau\in H^{\mathbb{S}}_{\dive}(Q_T)$. 

In view of the Signorini boundary conditions, the right-hand side of (\ref{eq: 33.1}) vanishes if $v=u$ and $\tau=\nabla u$.

\section{Examples}

In this section, we consider several examples that demonstrate how the estimate (\ref{eq:mainest}) works in practice. Namely, we consider the slightly modified  model  problem taken from \cite{NSV2000}, where the exact solution is known.
Due to this it is possible to examine the efficiency of error estimates for different approximate solutions $v$ and different
vector-functions $\tau$. 

\begin{figure}[htbp]
\centering
\includegraphics[width=0.44\textwidth]{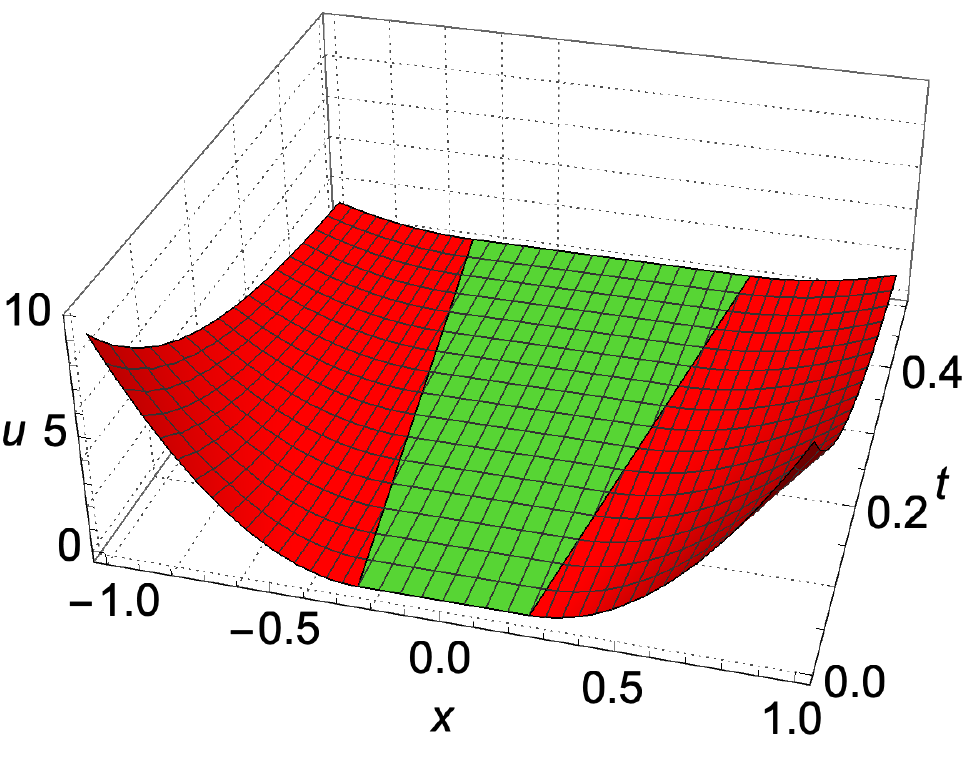}\qquad
\includegraphics[width=0.44\textwidth]{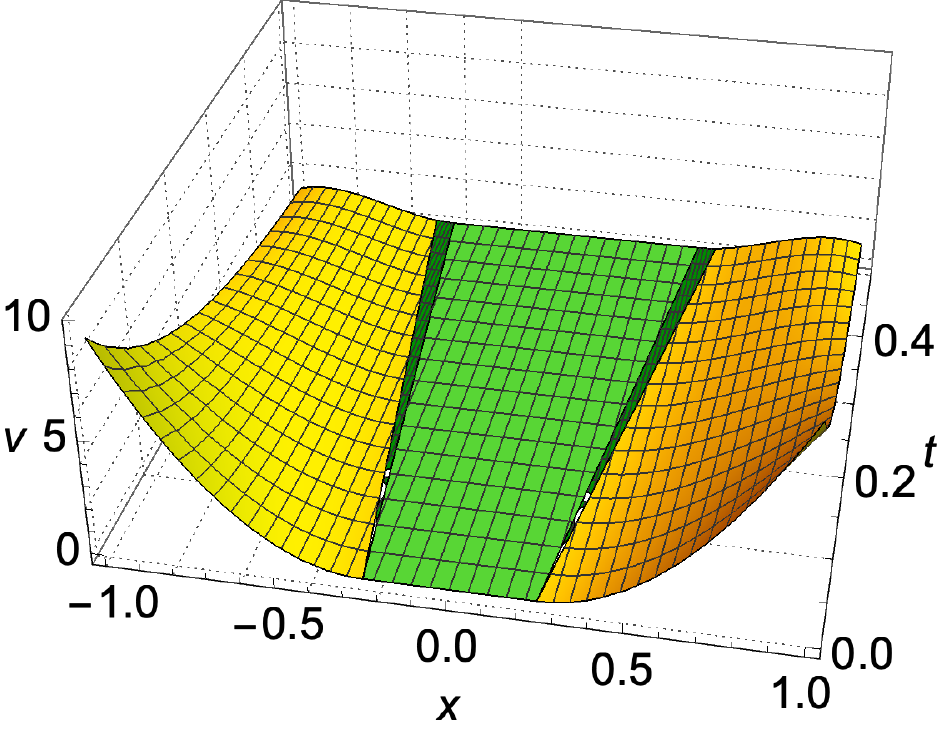}
\caption{The exact  solution $u (x,t)$ (left) and approximate solution $v_{\epsilon}(x,t)$  for $\varepsilon=0.5$ (right).}
\label{fig:1d_solution}
\end{figure}

\begin{figure}[!h]
\centering
\includegraphics[width=0.38\textwidth]{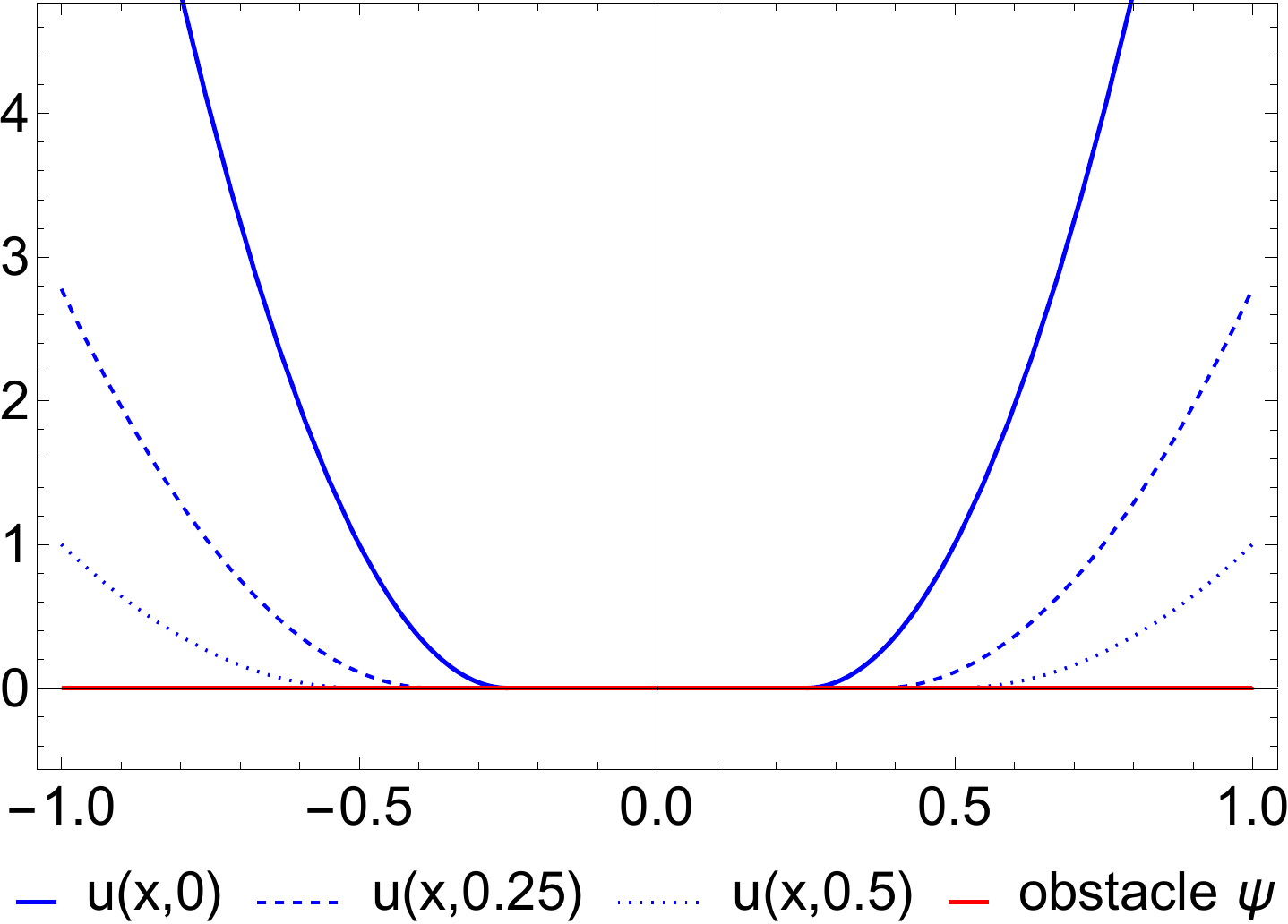}\quad
\includegraphics[width=0.49\textwidth]{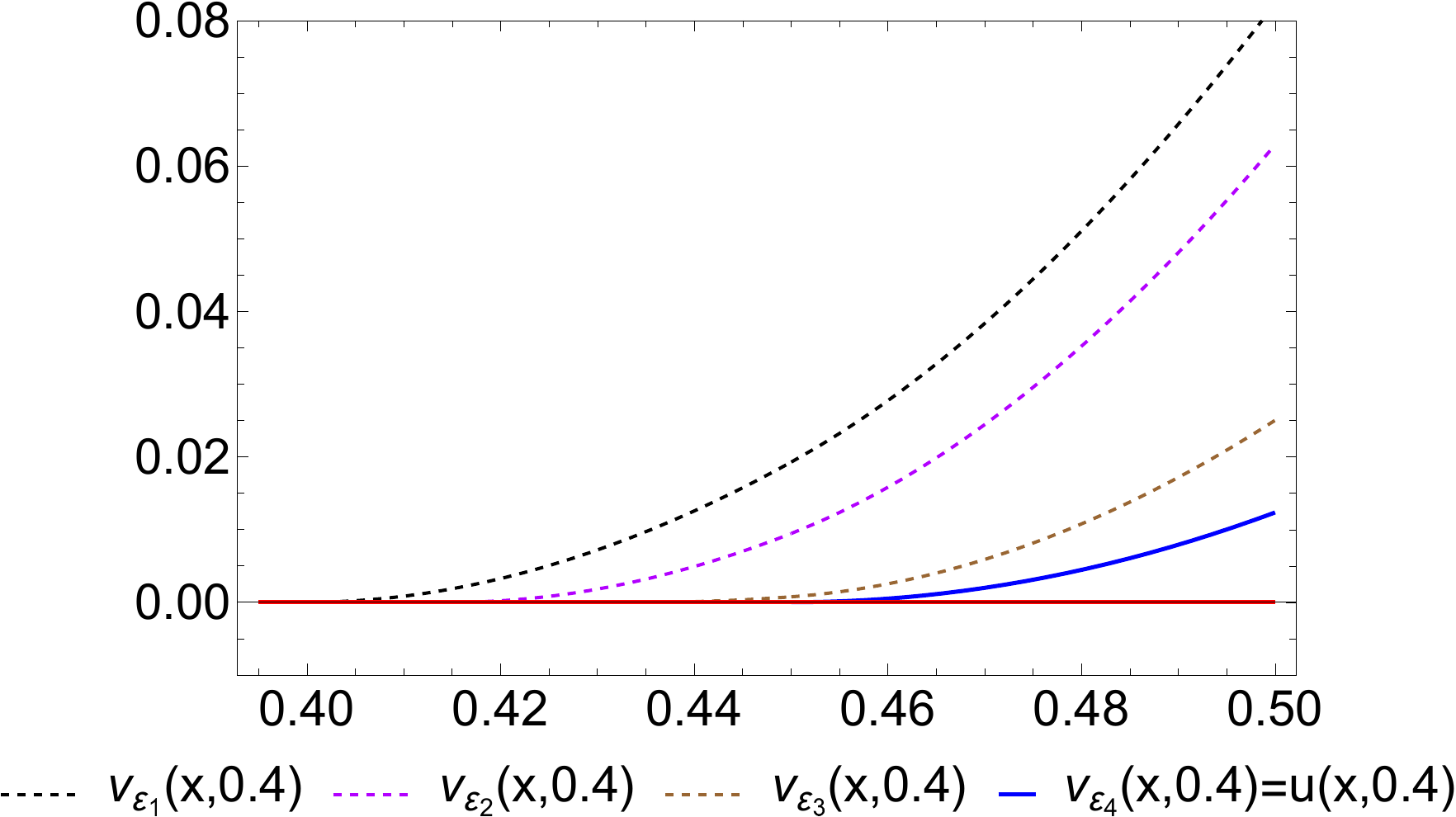}
\caption{The exact  1D-solution $u (\cdot,t)$  at times $t=0$, $0.25$, and $0.5$ (left)  and fragments  of the functions\\ $v_{\varepsilon} (\cdot,0.4)$ for $\varepsilon_1=0.5$, $\varepsilon_2=0.35$,   $\varepsilon_3=0.15$, and $\varepsilon_4=0$ (right).
}
\label{fig:v_e_sol}
\end{figure}

Let $\Omega=(-1.0,1.0)$,   $T=0.5$,  and $\phi \equiv 0$.
On the lateral surface $\partial''Q_T=\partial\Omega \times (0,T)$ we impose the boundary conditions
\begin{equation*} 
u(-1,t)=u(1,t)=\frac{9-12t+4t^2}{(2t+1)^2}, \qquad \quad \forall t\in (0,0.5),
\end{equation*} 
and  set 
\begin{equation*} 
u(x,0)=\left\{
\begin{aligned}
&16 x^2-8|x|+1, && \frac{1}{4}<|x|<1,\\
&0, && |x| \leqslant \frac{1}{4}.
\end{aligned}
\right.
\end{equation*}
If
\be
f(x,t)&=\left\{
\begin{array}{cc}
-\frac{16}{(2t+1)^2}\left(\frac{4}{2t+1}x^2-|x|+2\right), & (x,t) \in N:=\{4|x|>2t+1\},\\
0, & (x,t) \in \Lambda:=\{4|x| \leqslant 2t+1\}
\end{array}
\right.
\ee
then the exact solution is defined by the relation
\be
\begin{aligned}
u(x,t)&=\left\{
\begin{aligned}
&\frac{16}{(2t+1)^2}x^2-\frac{8}{2t+1}|x|+1, && (x,t) \in N\\
&0, && (x,t) \in \Lambda,
\end{aligned}
\right. 
\end{aligned}
\ee
The
function $u$ is depicted in Fig.~\ref{fig:1d_solution} (left) and in Fig.~\ref{fig:v_e_sol} (left) at  $t=0$, $0.25$, and $0.5$. The  coincidence set $Q_T^{\phi}(u)$  in 
Fig.~\ref{fig:1d_solution} is highlighted in light green.

Consider a set of approximations
$$
v_{\varepsilon}(x,t)
=\left\{
\begin{aligned}
&u(x,t)+100\varepsilon t (1-|x|)\bigg(x-\sgn{x}\frac{(2-\varepsilon)t+1}{4}\bigg)^2, && (x,t) \in N_{\varepsilon},\\
&0, &&(x,t) \in \Lambda_{\varepsilon},
\end{aligned}
\right.
$$
where the sets $N_{\varepsilon}$ and $\Lambda_{\varepsilon}$ are defined as follows:
\begin{align*}
N_{\varepsilon}&:=\{(x,t)\in Q_T: 4|x| > (2-\varepsilon)t+1\},\\
\Lambda_{\varepsilon}&:= \{(x,t)\in Q_T: 4|x| \leqslant (2-\varepsilon)t+1\}.
\end{align*}
Approximations depend on the parameter
 $\varepsilon\in [0, 1/2]$. The function $v_{\varepsilon}$ for $\varepsilon=0.5$ is depicted in Fig.~\ref{fig:1d_solution} (right). The coincidence set $Q_T^{\psi}(v)$ is marked in light green, while the values of $v_{\varepsilon}$ corresponding for $(x,t)\in \Lambda \setminus \Lambda_{\varepsilon}$ are highlighted in dark green.

For any $\varepsilon$, the function
$v_{\varepsilon}$ belongs to the set  ${\mathbb{K}}$ and $\Lambda_{\varepsilon} \subset \Lambda$. If $\varepsilon\rightarrow 0$,
then $\Lambda_\varepsilon$ tends to
$\Lambda=\{(x,t) \mid u=\psi\}$ and  $v_\varepsilon$ tends to $u$ (see, Fig.~\ref{fig:v_e_sol} (right)).

\begin{figure}[htbp]
\centering
\includegraphics[width=0.93\textwidth]{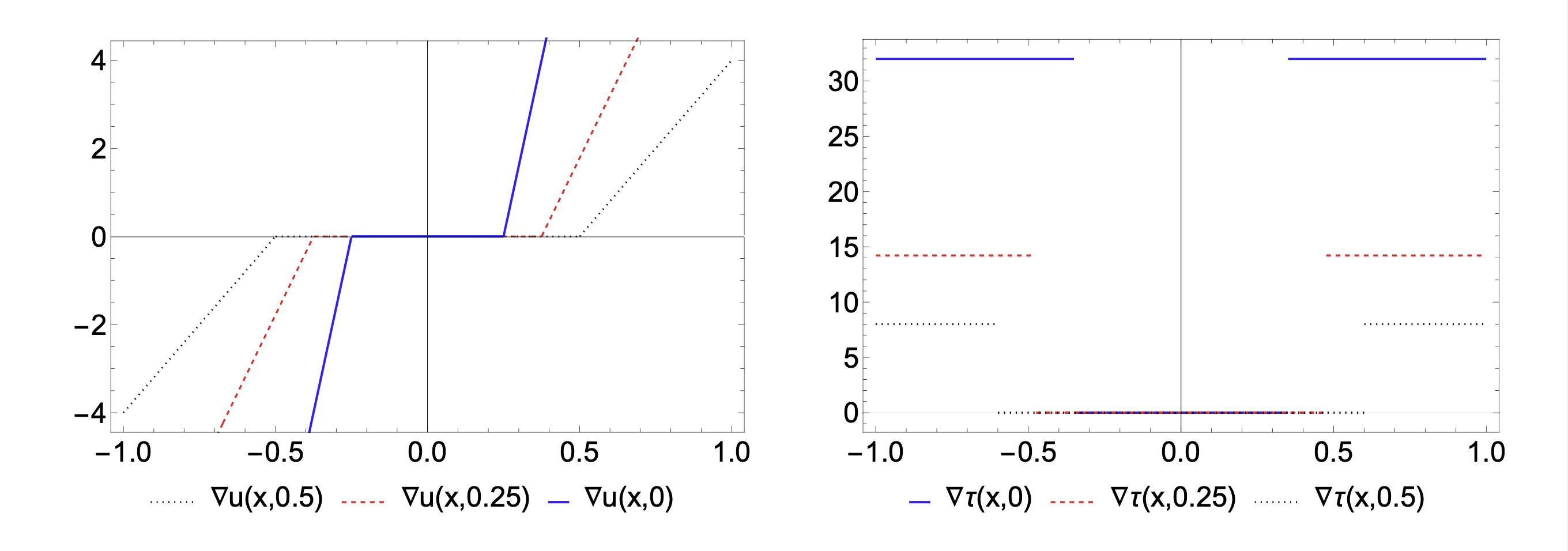}
\caption{The functions $\tau (\cdot,t)=\nabla u(\cdot, t)$ (left) and $\tau _x(\cdot,t)$ (right) at times $t=0, 0.25$, and $0.5$.}
\label{fig:tau_bild}
\end{figure}

\vspace{0.2cm}
First, we set
$$
\tau  (x,t)= \left\{
\begin{aligned}
&\frac{32x}{(2t+1)^2}-\sgn{x}\cdot\frac{8}{2t+1}, && (x,t) \in N,\\
&0, && (x,t) \in \Lambda.
\end{aligned}
\right.
$$
So defined  $\tau$ belongs to  $L^2(-1,1)$ ( Fig.~\ref{fig:tau_bild}), and, consequently, $\tau \in H_{\dive}(Q_T)$. 


In our case, $\Omega=(-1,1)$ and $C_F=2/\pi$. We verify the validity of the estimate (\ref{eq:mainest}) for $v=v_{\epsilon}$ with different $\varepsilon$. Table~\ref{tb:table1} presents the results related to the components of (\ref{eq:mainest}) for $\varepsilon= 0.05j$, $j=10, 7, 5, 3,$ and $0$. It shows how different terms of the error measure
and error majorant decrease as $\varepsilon \to 0$.
\begin{table}[!h] 
\centering
\begin{tabular}{c|ccccc}
$\varepsilon $  & $\|e(\cdot,0.5)\|^2_{\Omega}$ & $\|\nabla e\|^2_{Q_T}$ & $\|e(\cdot,0)\|^2_{\Omega}$ & $\|\tau-\nabla v_{\varepsilon} \|_{Q_T}$  & $\|\mathcal{F}_f(v_{\varepsilon},\tau )\|_{Q_T}$ \\
\midrule
0.50  & $21.21\cdot 10^{-2}$  & $2.42$  & 0  & 1.56  & 0.87\\
0.35  & $8.20\cdot 10^{-2}$  & $1.07$  & 0  & 1.03  & 0.59\\
0.25  & $3.55\cdot 10^{-2}$  & $0.51$  & 0  & 0.71  & 0.41\\
0.15  & $1.08\cdot 10^{-2}$  & $0.17$  & 0  & 0.41  & 0.24 \\
0.05  & $1.02\cdot 10^{-3}$  & $1.75\cdot 10^{-2}$  & 0  & 0.13  & $7.87\cdot 10^{-2}$\\
0.00  & 0  & 0  & 0  & 0 &  0\\
\end{tabular}
\caption{Components of the estimate (\ref{eq:mainest}) for $v=v_{\epsilon}$ and $\tau =\nabla u$ for different $\varepsilon$.}
\label{tb:table1}
\end{table}

Table~\ref{tb:table2} presents the results in the integral form.
It compares  exact errors (l.h.s. of (\ref{eq:mainest})) and  error majorants (r.h.s. of (\ref{eq:mainest})) computed for $\alpha=1$, and $2$
together with the corresponding efficiency indices 
$$
1 \leqslant I_{\rm{eff}}=\sqrt{\frac{\rm{r.h.s.\ of} (\ref{eq:mainest})}{\rm{l.h.s.\ of (\ref{eq:mainest})}}}.
$$
We see that the results are very good, what one may await because
in these tests  the
function $\tau$ coincides with the exact gradient $\nabla u$.
\vspace{0.5cm}
\begin{table}[!h] 
\centering
\begin{tabular}{c|cccc||cccc}
{}  &\multicolumn{3}{c}{$\alpha=1$} &{} &\multicolumn{3}{c}{$\alpha=2$}  \\
\midrule
$\varepsilon $ & l.h.s. of (\ref{eq:mainest}) & r.h.s. of (\ref{eq:mainest})  &$I_{\rm{eff}}$ & \, & \ l.h.s. of (\ref{eq:mainest}) & r.h.s. of (\ref{eq:mainest})  &$I_{\rm{eff}}$   \\
\midrule
0.50 & 2.63  & 4.47  & 1.304 &\ & 3.84 &  8.94 & 1.526 \\ 
0.35  & 1.15  & 1.98  & 1.312 &&  1.69 & 3.95 & 1.529 \\
0.25 & 0.55  & 0.94  & 1.307 && 0.80 & 1.89 & 1.537 \\
0.15  & 0.18  & 0.32  & 1.333 &&
0.27 & 0.63 & 1.528 \\
0.05  & $1.85\cdot 10^{-2}$  & $3.24\cdot 10^{-2}$  & 1.323 && $2.73 \cdot 10^{-2}$& $6.49 \cdot 10^{-2}$ & 1.542  \\
0  & 0  & 0  & \ && 0 & 0 && \\
\end{tabular}
\caption{Estimate (\ref{eq:mainest}) for $v=v_{\varepsilon}$ and $\tau =\nabla u$ with $\alpha=1$ and $\alpha=2$ for different $\varepsilon$.}
\label{tb:table2}
\end{table}


\vspace{-0.5cm}
\begin{figure}[!h]
\centering
\includegraphics[width=0.95\textwidth]{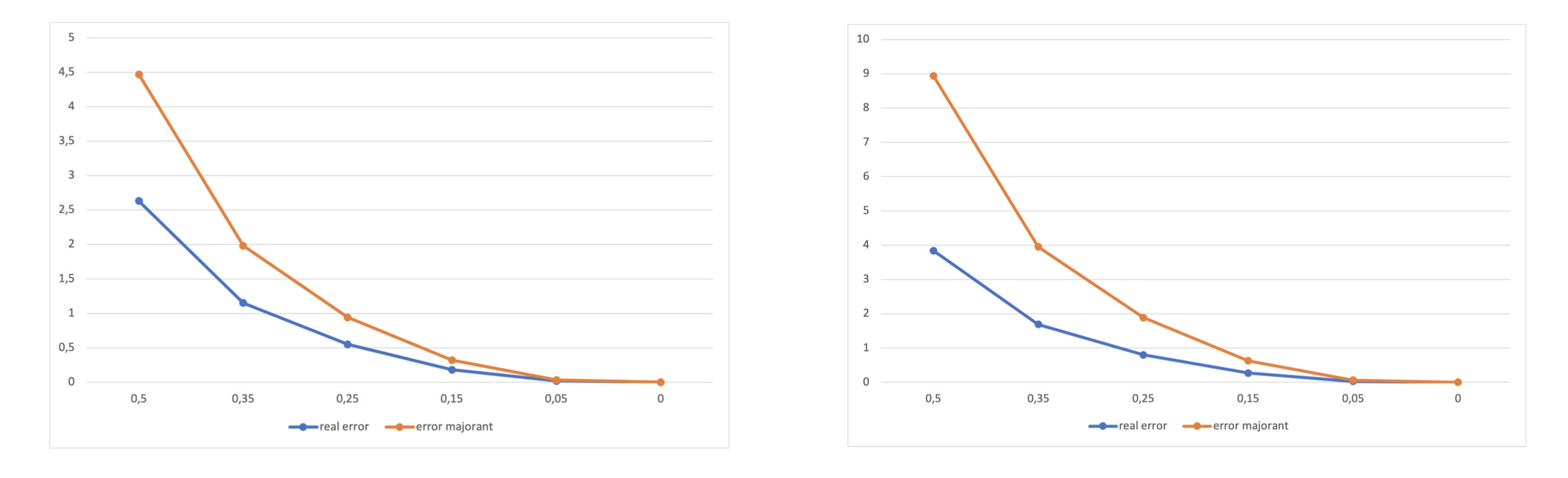}
\caption{The real error  (l.h.s. of (\ref{eq:mainest}))  vs the error majorant  (r.h.s. of (\ref{eq:mainest})) computed  for $\alpha=1$ (left)\\ and $\alpha=2$ (right).}
\label{fig:error-vs-majorant}
\end{figure}
Graphically, these results are depicted
in Fig.~\ref{fig:error-vs-majorant}.

\begin{figure}[!h]
\centering
\includegraphics[scale=1.4]{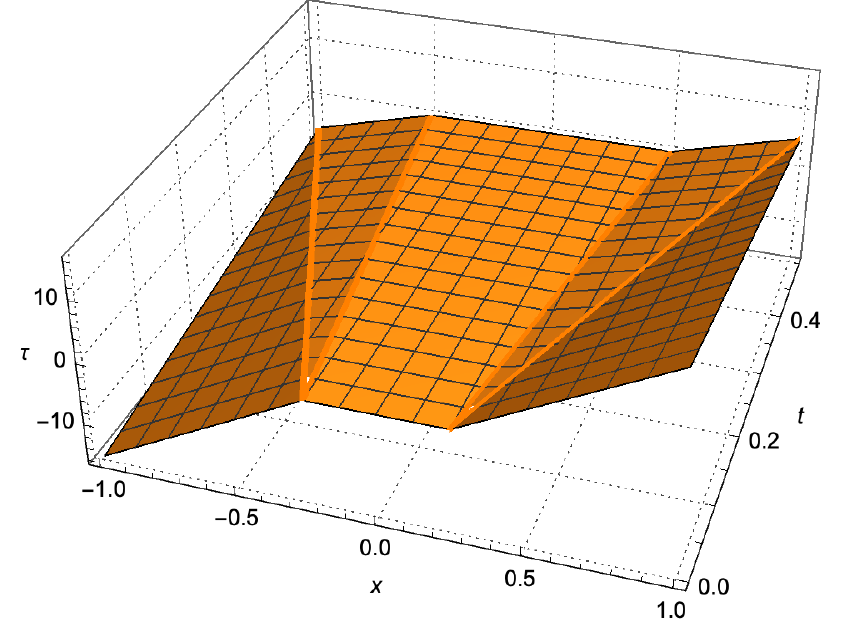}
\caption{The graph of $\tau_{\delta}$  for $\delta=0.5$, $\xi=16.07$, and $\eta=5.62$.}
\label{fig:T_delta}
\end{figure}

Further, we investigate the question on ability
to get close results using relatively simple approximations
of $\tau$ that contain only constants and terms proportional to $t$ and $x$. Certainly, so simple functions may be useful only for a small
time interval (e.g., for time incremental methods that operate
with small intervals). Therefore, we set the time parameter $\delta \in (0,0.5]$ and define the space--time domain $Q_{\delta}:=\Omega\times (0,\delta)$. Consider   a sequence of  functions $\tau=\tau_{\delta}(x,t)$ defined by the relation
\begin{equation} \label{eq:tau_d}
\tau_{\delta}(x,t)=
\left\{
\begin{aligned}
&0, && \text{in}\  Q_{\delta}\cap\left\{|x| \leqslant \frac{1+2t}{4}\right\}\\
&\frac{4\eta}{3-2\delta}\left(x- \sgn{x}\cdot\frac{1+2t}{4}\right),
&& \text{in}\ Q_{\delta}\cap\left\{\frac{1+2t}{4}<|x|\leqslant \frac{\delta +3t}{4\delta}\right\},\\
&\frac{4\xi}{3}\left(x-\frac{\sgn{x}}{4}\right)+\sgn{x}(\eta-\xi)\frac{t}{\delta}, && \text{in}\ Q_{\delta}\cap\left\{\frac{\delta+3t}{4\delta}<|x|\leqslant 1\right\}.
\end{aligned}
\right.
\end{equation}
The values of  coefficients $\xi=\xi(\delta)$ and $\eta=\eta (\delta)$  in formula (\ref{eq:tau_d}) are determined by minimising the corresponding right-hand sides of  (3.3). One can easily check that for all $\delta$ the functions $\tau_{\delta}$ are  continuous (see, for example, Fig.~\ref{fig:T_delta}).  Therefore, $\tau_{\delta}$ again belongs to the set $H_{\dive}(Q_{\delta})$. 

\vspace{0.2cm}
Table~\ref{tb:table3}  reports on values of the real errors and the error majorants included in (\ref{eq:mainest}) for few approximations $v_{\epsilon}$, where $\epsilon=0.2$.
\begin{table}[!h] 
\centering
\begin{tabular}{c|cccccccccc}
$\delta$ & \ & $\xi (\delta)$& \ & $\eta (\delta)$ &\ & l.h.s. of (\ref{eq:mainest}) &\ & r.h.s. of (\ref{eq:mainest}) &\ & $I_{\rm{eff}}$ \\
\midrule
0.5 & \ & 16.07 &\ & 5.62 &\ & 0.33 &\ & 19.89&\ & 7.722\\
0.3 & \ & 18.68 &\ & 9.61 &\ & 0.13 &\ & 8.27 &\ & 7.829\\
0.2 & \ & 20.3 &\ & 12.78 &\ & 0.06 &\ & 3.64 &\ & 7.857\\
0.1 & \ & 22.18 & \ & 17.33 & \ & 0.01 &\ & 0.72 &\ & 8.487\\
\end{tabular}
\caption{Estimate (\ref{eq:mainest}) for $v=v_{\epsilon}$, $\varepsilon=0.2$, $\tau =\tau_{\delta}$, $\alpha=1$, and different $\delta$.}
\label{tb:table3}
\end{table}

Table~\ref{tb:table4} collects values of left- and right-hand sides of (\ref{eq:mainest}) for $\tau=\tau_{\delta}$ and  for time-incremental approximations $v=w_{\delta}(x,t):=u(x,0)+\frac{u(x,\delta)-u(x,0)}{\delta}t$
for $(x,t)\in Q_{\delta}$. Here the estimates are coarser, what is
not surprising because rather coarse approximations are estimated
using a simple function $\tau_\delta$. Anyway, even in this case, the majorant gives a presentation on the actual value of the error.

\begin{table}[!h] 
\centering
\begin{tabular}{c|cccccccccc}
$\delta$ & \ & $\xi (\delta)$& \ & $\eta (\delta)$ &\ & l.h.s. of (\ref{eq:mainest}) &\ & r.h.s. of (\ref{eq:mainest}) &\ & $I_{\rm{eff}}$ \\
\midrule
0.5 & \ & 18.08 &\ & 7.22 &\ & 4.54 &\ & 27.41&\ & 2.456 \\
0.3 & \ & 21.54 &\ & 10.23 &\ & 0.89 &\ & 22.43 &\ & 5.024 \\
0.2 & \ & 22.38 &\ & 12.92 &\ & 0.20 &\ & 9.86 &\ & 6.963 \\
\end{tabular}
\caption{Estimate (\ref{eq:mainest}) for $v=w_{\delta}$, $\tau =\tau_{\delta}$,  $\alpha=1$, and different $\delta$.}
\label{tb:table4}
\end{table}

\vskip-0.2cm
\begin{figure}[!h]
\centering
\includegraphics[scale=1.4]{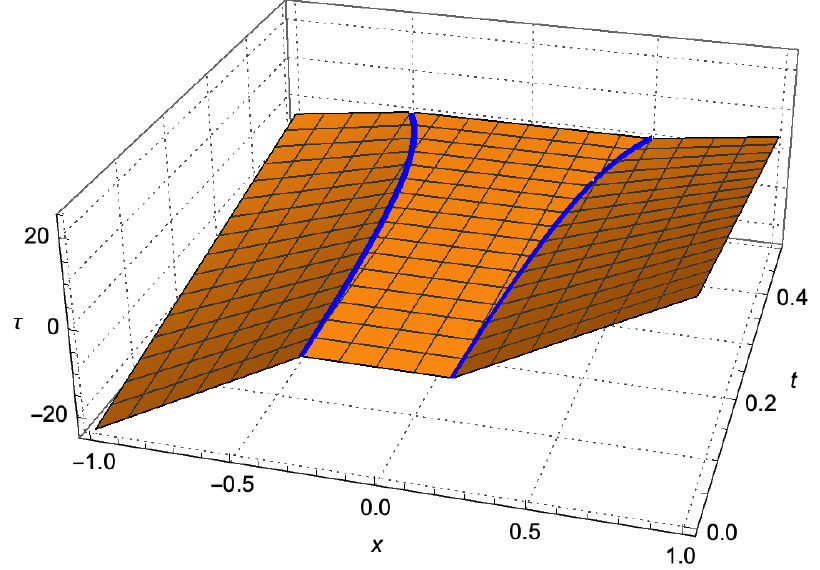}
\caption{The graph of $\hat{\tau}_{\delta}$  for $\delta=0.5$, $\xi=24$, and $\eta=5.62$.}
\label{fig:T_polylinear}
\end{figure}

In the final series of tests, we used
\begin{equation*} 
\hat{\tau}_{\delta}(x,t)=
\left\{
\begin{aligned}
&32\left(x- \sgn{x}\frac{1+2t}{4}\right)-\frac{128(1+\delta)t}{(1+2\delta )^2}\left(x-\sgn{x}\frac{1+2\delta}{4}\right),
&& (x,t)\in N_{\delta},\\
&0, && (x,t)\in \Lambda_{\delta}\\
\end{aligned}
\right.
\end{equation*}
where $\delta \in (0,0.5]$ and sets are defined as follows:
\begin{align*}
N_{\delta}&:=Q_{\delta}\cap\left\{\frac{(1+2\delta)(1+2\delta-2t)}{4(1+4\delta^2-4\delta (t-1)-4t)}<|x|\leqslant 1\right\},\\
\Lambda_{\delta}&:=Q_{\delta} \cap \left\{|x| \leqslant \frac{(1+2\delta)(1+2\delta-2t)}{4(1+4\delta^2-4\delta (t-1)-4t)}\right\}.
\end{align*}

For all $\delta$, the functions $\hat{\tau}_{\delta}$ are continuous (see, e.g. Fig.~\ref{fig:T_polylinear}), and $\left(\hat{\tau}_{\delta}\right)_x \in L^2(-1,1)$. Thus, the condition $\hat{\tau}_{\delta} \in H_{\dive}(Q_{\delta})$ is fulfilled as well.

Table~\ref{tb:table5} shows  values of the exact errors and respective majorants from (\ref{eq:mainest}) computed for two types of the approximate solutions: $v=v_{\epsilon}$ and $v=w_{\delta}$. Comparing the results 
with Tables~\ref{tb:table3}--\ref{tb:table5}, we see that using $\hat{\tau}_{\delta}$ instead of $\tau=\tau_{\delta}$
improves the estimates, especially for small values of $\delta$.

\vspace{-0.2cm}
\begin{table}[!h] 
\centering
\begin{tabular}{c|cccc||cccc}
{}  &\multicolumn{3}{c}{$v=v_{\varepsilon}, \quad \varepsilon=0.2$} &{} &\multicolumn{3}{c}{$v=w_{\delta}$}  \\
\midrule
$\delta $  & l.h.s. of (\ref{eq:mainest})& r.h.s. of (\ref{eq:mainest})   &$I_{\rm{eff}}$ & \quad & \ l.h.s. of (\ref{eq:mainest}) & r.h.s. of (\ref{eq:mainest}) &  $I_{\rm{eff}}$   \\
\midrule
0.5  &   &   &  &\ & 4.54 &  24.68 & 2.331 \\ 
0.3  & 0.13  & 8.43  & 7.901 &&  0.89 & 9.13 & 3.021 \\
0.2  & 0.06  & 1.58  & 5.179 && 0.20 & 3.93 & 1.983 \\
0.1  & 0.01  & 0.34  & 5.860 &&
 &&  &  \\
\end{tabular}
\caption{Estimate (\ref{eq:mainest}) for  $\tau =\hat{\tau}_{\delta}$ with $\alpha=1$ and  different $\delta$.}
\label{tb:table5}
\end{table}


 The first author was supported by the Deutsche Forschungsgemeinschaft (DFG, German Research Foundation) - Project-ID AP 252/3-1.
\bibliographystyle{alpha}
\bibliography{Bibliography(par-obstacle)}

\end{document}